\documentclass[12pt]{article}
\usepackage{amsmath,amssymb,amsthm}
\numberwithin{equation}{section}
\usepackage{hyperref}
\usepackage{units}
\usepackage{color}
\usepackage[T1]{fontenc}
\usepackage[utf8]{inputenc}
\usepackage{authblk}
\usepackage{bm}
\usepackage{graphicx}

\usepackage{datetime}

\usepackage[toc,page]{appendix}

\newcommand{\Z}{{\mathbb Z}}

\newtheorem{thm}{Theorem}

\newtheorem{lemma}{Lemma}

\title{Weighted sums of some second-order sequences}

\author[]{Kunle Adegoke \thanks{adegoke00@gmail.com}}

\affil{Department of Physics and Engineering Physics, \mbox{Obafemi Awolowo University}, 220005 Ile-Ife, Nigeria}

\begin{document}

\date{}

\maketitle

\begin{abstract}
\noindent We derive weighted summation identities involving the second order recurrence sequence $\{w_n\} =\{ w_n(a,b; p, q)\}$ defined by $w_0  = a,\,w_1  = b;\,w_n  = pw_{n - 1}  - qw_{n - 2}\, (n \ge 2)$, where  $a$, $b$, $p$ and $q$ are arbitrary complex numbers, with $p\ne 0$ and $q\ne 0$. 
\end{abstract}

\section{Introduction}
Horadam~\cite{horadam65} wrote a fascinating paper in which he established the basic arithmetical properties of his generalized Fibonacci sequence $\{w_n\} = \{w_n(a,b; p, q)\}$ defined by
\begin{equation}
w_0  = a,\,w_1  = b;\,w_n  = pw_{n - 1}  - qw_{n - 2}\, (n \ge 2)\,,
\end{equation}
where  $a$, $b$, $p$ and $q$ are arbitrary complex numbers, with $p\ne 0$ and $q\ne 0$. 
Some well studied particular cases of $\{w_n\}$ are the sequences $\{u_n\}$, $\{v_n\}$, $\{G_n\}$, $\{P_n\}$ and $\{J_n\}$ given by:
\begin{equation}
w_n (1,p;p,q) = u_n (p,q)\,,
\end{equation}
\begin{equation}
w_n (2,p;p,q) = v_n (p,q)\,,
\end{equation}
\begin{equation}
w_n (a,b;1,-1) = G_n (a,b)\,,
\end{equation}
\begin{equation}
w_n (0,1;2,-1) = P_n
\end{equation}
and
\begin{equation}
w_n (0,1;1,-2) = J_n\,.
\end{equation}
Note that \mbox{$u_n(1,-1)=F_{n+1}$} and $v_n(1,-1)=L_n$, where $F_n=G_n(0,1)$ and $L_n=G_n(2,1)$ are the classic Fibonacci numbers and Lucas numbers, respectively. $P_n$ and $J_n$ are the Pell numbers and Jacobsthal numbers, respectively. Note also that $u_n(2,-1)=P_{n+1}$ and \mbox{$u_n(1,-2)=J_{n+1}$}. The sequence $\{G_n\}$ was introduced by Horadam~\cite{horadam61} in 1961, (under the notation $\{H_n\}$).

\bigskip

Extension of the definition of $w_n$ to negative subscripts is provided by writing the recurrence relation as $w_{-n}=(pw_{-n+1}-w_{-n+2})/q$.
Horadam~\cite{horadam65} showed that:
\begin{equation}
u_{-n}=-q^{-n+1}u_{n-2}\,,
\end{equation}
\begin{equation}
v_{-n}=q^nv_n
\end{equation}
and
\begin{equation}
w_{ - n}  = \frac{{au_n  - bu_{n - 1} }}{{au_n  + (b - pa)u_{n - 1} }}w_n\,.
\end{equation}
Our main goal in this paper is to derive weighted summation identities involving the numbers $w_n$. For example, we shall derive (Theorem~\ref{thm.g1ihfq5}) the following weighted binomial sum:
\[
(-qu_{r - 1})^k\sum_{j = 0}^k {\binom kj\left( {-\frac{{u_r }}{{qu_{r - 1} }}} \right)^j w_{m - k(r + 1) + j} }  = w_m\,,
\]
which generalizes Horadam's result (\cite{horadam65}, equation 3.19):
\[
( - q)^n\sum_{j = 0}^n {\binom nj\left( { - \frac{p}{q}} \right)^j w_j }  =  w_{2n}\,,
\]
the latter identity being an evaluation of the former at $m=2n$, $k=n$, $r=1$.

\bigskip

As another example, it is known (first identity of Cor.~15, \cite{stanica03}) that
\[
\sum_{j = 0}^n {( - 1)^j \binom nj F_j }  =  - F_n\,, 
\]
but this can be generalized to:
\[
\sum_{j = 0}^k { (-1)^j  \binom kj\frac{G_{rj}}{F_{r + 1}^j} }  = \left( {\frac{{F_r }}{{F_{r + 1} }}} \right)^k \frac{{G_0 F_{k + 1}  - G_1 F_k }}{{G_0 F_{k - 1}  + G_1 F_k }}G_k\,,
\]
which is itself a special case of a more general result (see Theorem~\ref{thm.peyb26i}):
\[
\sum_{j = 0}^k {(-1)^j\binom kj \frac{w_{rj}}{u_r^j} }  = \left( {\frac{{u_{r - 1} }}{{u_r }}} \right)^k \frac{{au_k  - bu_{k - 1} }}{{au_k  + (b - pa)u_{k - 1} }}w_k\,.
\]
As an example of non-binomial sums derived in this paper we mention (see Theorem~\ref{thm.m4yo6a8}):
\[
u_{r-1} \sum_{j = 0}^k {\frac{{w_{rj} }}{(-qu_{r - 2})^j}}  =  \frac{{w_{kr + r - 1} }}{(-qu_{r - 2})^k} + (ap - b)u_{r - 2}\,,
\]
of which a special case is
\[
F_r \sum_{j = 0}^k {\frac{{G_{rj} }}{{F_{r-1}^j }}}  = \frac{{G_{kr + r - 1} }}{{F_{r-1}^k }}- F_{r-1} (G_1  - G_0 )\,.
\]
Another example in this category is (see Theorem~\ref{thm.ybopnqn}):
\[
q^{n-r} eu_{r - 1} \sum_{j = 0}^k {\frac{{u_{rj} }}{{(w_n /w_{n - r} )^j }}}  = \frac{{w_{n + kr + 1} w_{n - r} }}{{(w_n /w_{n - r} )^k }} - w_n w_{n - r + 1}\,,
\]
giving rise to the following results for the $\{G_m\}$, $\{P_m\}$ and $\{J_m\}$ sequences:
\[
( - 1)^{n - r} (G_0 G_1  + G_0^2  - G_1^2 )F_r \sum_{j = 0}^k {\frac{{F_{rj + 1} }}{{(G_n /G_{n - r} )^j }}}  = \frac{{G_{n + kr + 1} G_{n - r} }}{{(G_n /G_{n - r} )^k }} - G_n G_{n - r + 1}\,,
\]
\[
( - 1)^{n - r - 1} P_r \sum_{j = 0}^k {\frac{{P_{rj + 1} }}{{(P_n /P_{n - r} )^j }}}  = \frac{{P_{n + kr + 1} P_{n - r} }}{{(P_n /P_{n - r} )^k }} - P_n P_{n - r + 1}
\]
and
\[
( - 1)^{n - r - 1}2^{n-r} J_r \sum_{j = 0}^k {\frac{{J_{rj + 1} }}{{(J_n /J_{n - r} )^j }}}  = \frac{{J_{n + kr + 1} J_{n - r} }}{{(J_n /J_{n - r} )^k }} - J_n J_{n - r + 1}\,.
\]
We require the following identities, derived in~\cite{horadam65}:
\begin{equation}\label{eq.fuxige6}
w_{m+r}=u_rw_m-qu_{r-1}w_{m-1}\,,
\end{equation}
\begin{equation}\label{eq.w7u7hr6}
v_rw_m=w_{m+r}+q^rw_{m-r}
\end{equation}
and
\begin{equation}\label{eq.vx6b1t4}
w_{n - r} w_{m + n + r}  = w_n w_{m + n}  + q^{n - r} eu_{r - 1} u_{m + r - 1}\,,
\end{equation}
where $e=pab-qa^2-b^2$.
\section{Weighted sums}
\begin{lemma}\label{lem.u4bqbkc}
Let $\{X_m\}$ and $\{Y_m\}$ be any two sequences such that $X_m$ and $Y_m$, $m\in\Z$, are connected by a second order recurrence relation $X_m=xX_{m-\alpha}+yY_{m-\beta}$, where $x$ and $y$ are arbitrary non-vanishing complex functions, not dependent on $m$, and $\alpha$ and $\beta$ are integers. Then,
\[
y\sum_{j = 0}^k {\frac{{Y_{m - k\alpha  - \beta  + \alpha j} }}{{x^j }}}  = \frac{{X_m }}{{x^k }} - xX_{m - (k + 1)\alpha }\,, 
\]
for $k$ a non-negative integer.
\end{lemma}
In particular,
\begin{equation}
y\sum_{j = 0}^k {\frac{{Y_{\alpha j} }}{{x^j }}}  = \frac{{X_{k\alpha  + \beta } }}{{x^k }} - xX_{\beta  - \alpha }\,.
\end{equation}
\begin{proof}
The proof shall be by induction on $k$. Consider the proposition $\mathcal{P}_k$,
\[
\mathcal{P}_k:\left(y\sum_{j = 0}^k {\frac{{Y_{m - k\alpha  - \beta  + \alpha j} }}{{x^j }}}  = \frac{{X_m }}{{x^k }} - xX_{m - (k + 1)\alpha }\right)\,; 
\]
with respect to the relation $X_m=xX_{m-\alpha}+yY_{m-\beta}$.
Clearly, $\mathcal{P}_0$ is true. Assume that $\mathcal{P}_n$ is true for a certain positive integer $n$. We want to prove that $\mathcal{P}_n\Rightarrow\mathcal{P}_{n+1}$.
Now,
\[
\mathcal{P}_n:\left(f(n)  = \frac{{X_m }}{{x^n }} - xX_{m - (n + 1)\alpha }\right)\,; 
\]

where
\[
f(n) = y\sum_{j = 0}^n {\frac{{Y_{m - n\alpha  - \beta  + j\alpha } }}{{x^j }}}\,. 
\]
We have
\[
f(n + 1) = y\sum_{j = 0}^{n + 1} {\frac{{Y_{m - n\alpha  - \alpha  - \beta  + j\alpha } }}{{x^j }}}\qquad\qquad\quad
\]
\[
\begin{split}
&= y\sum_{j =  - 1}^n {\frac{{Y_{m - n\alpha  - \alpha  - \beta  + j\alpha  + \alpha } }}{{x^{j + 1} }}}\\
&= \frac{y}{x}\sum_{j =  - 1}^n {\frac{{Y_{m - n\alpha  - \beta  + j\alpha} }}{{x^j }}}
\end{split}
\]
\[
\begin{split}
&\qquad\qquad\qquad\quad= \frac{y}{x}\left( {xY_{m - n\alpha  - \beta  - \alpha }  + \sum_{j = 0}^n {\frac{{Y_{m - n\alpha  - \beta  + j\alpha } }}{{x^j }}} } \right)\\
&\qquad\qquad\qquad\quad= yY_{m - n\alpha  - \alpha  - \beta }  + \frac{1}{x}\left(y\sum_{j = 0}^n {\frac{{Y_{m - n\alpha  - \beta  + j\alpha } }}{{x^j }}}\right)
\end{split}
\]
\qquad\qquad\qquad\qquad\qquad(invoking the induction hypothesis $\mathcal{P}_n$)
\[
\begin{split}
&\qquad\qquad\qquad= yY_{m - n\alpha  - \alpha  - \beta }  + \frac{1}{x}\left( {\frac{{X_m }}{{x^n }} - xX_{m - n\alpha  - \alpha } } \right)\\
&\qquad\qquad\qquad= \frac{{X_m }}{{x^{n + 1} }} - \left( {X_{m - n\alpha  - \alpha }  - yY_{m - n\alpha  - \alpha  - \beta } } \right)\,.
\end{split}
\]
Since $X_{m - n\alpha  - \alpha }  - yY_{m - n\alpha  - \alpha  - \beta }=xX_{m - n\alpha  - \alpha - \alpha} $, we finally have
\[
f(n+1) = \frac{{X_m }}{{x^{n + 1} }} - xX_{m - (n + 2)\alpha }\,.
\]
Thus,
\[
\mathcal{P}_{n+1}:\left(f(n+1)  = \frac{{X_m }}{{x^{n+1} }} - xX_{m - (n + 1+1)\alpha }\right)\,; 
\]
i.~e.~$\mathcal{P}_n\Rightarrow\mathcal{P}_{n+1}$ and the induction is complete.
\end{proof}
Note that the identity of Lemma~\ref{lem.u4bqbkc} can also be written in the equivalent form:
\begin{equation}\label{eq.t347olg}
y\sum_{j = 0}^k {x^j Y_{m - \beta  - j\alpha } }  = X_m  - x^{k + 1} X_{m - (k + 1)\alpha }\,.
\end{equation}
In particular,
\begin{equation}
y\sum_{j = 0}^k {x^j Y_{ - j\alpha } }  = X_\beta   - x^{k + 1} X_{\beta  - (k + 1)\alpha }\,.
\end{equation}
\begin{thm}\label{thm.xvb2v42}
For integer $m$, non-negative integer $k$ and any integer $r$ for which $w_{r-1}\ne 0$, the following identity holds:
\[
\sum_{j = 0}^k {\left( {\frac{{w_r }}{{qw_{r - 1} }}} \right)^j w_{m + r - k + j} }  = \left( {\frac{{w_r }}{{qw_{r - 1} }}} \right)^k u_m w_r  - qu_{m - k - 1} w_{r - 1}\,.
\]

\end{thm}
In particular,
\begin{equation}
q^{r - 1} \sum_{j = 0}^k {\left( {\frac{{w_r }}{{qw_{r - 1} }}} \right)^j w_j }  = \left( {\frac{{w_r }}{{qw_{r - 1} }}} \right)^k q^{r - 1} u_{k - r} w_r  + u_{r - 1} w_{r - 1}\,.
\end{equation}
\begin{proof}
Interchange $m$ and $r$ in identity~\eqref{eq.fuxige6} and write the resulting identity as
\[
u_m  = \frac{{qw_{r - 1} }}{{w_r }}u_{m - 1}  + \frac{1}{{w_r }}w_{m + r}\,.
\]
Identify $X=u$, $Y=w$, $x=qw_{r-1}/w_r$, $y=1/w_r$, $\alpha=1$ and $\beta=-r$ and use these in Lemma~\ref{lem.u4bqbkc}.
\end{proof}
The Fibonacci, Lucas and Pell versions of Theorem~\ref{thm.xvb2v42} are, respectively,
\begin{equation}
\sum_{j = 0}^k {( - 1)^j \left( {\frac{{F_r }}{{F_{r - 1} }}} \right)^j F_{m + r - k + j} }  = ( - 1)^k \left( {\frac{{F_r }}{{F_{r - 1} }}} \right)^k F_{m + 1} F_r  + F_{m - k} F_{r - 1}\,,
\end{equation}
\begin{equation}
\sum_{j = 0}^k {( - 1)^j \left( {\frac{{L_r }}{{L_{r - 1} }}} \right)^j L_{m + r - k + j} }  = ( - 1)^k \left( {\frac{{L_r }}{{L_{r - 1} }}} \right)^k F_{m + 1} L_r  + F_{m - k} L_{r - 1}
\end{equation}
and
\begin{equation}
\sum_{j = 0}^k {( - 1)^j \left( {\frac{{P_r }}{{P_{r - 1} }}} \right)^j P_{m + r - k + j} }  = ( - 1)^k \left( {\frac{{P_r }}{{P_{r - 1} }}} \right)^k P_{m + 1} P_r  + P_{m - k} P_{r - 1}\,.
\end{equation}
In particular, we have
\begin{equation}
\sum_{j = 0}^k {( - 1)^j \left( {\frac{{F_r }}{{F_{r - 1} }}} \right)^j F_{r + j} }  = ( - 1)^k \left( {\frac{{F_r }}{{F_{r - 1} }}} \right)^k F_{k + 1} F_r\,,
\end{equation}
\begin{equation}
\sum_{j = 0}^k {( - 1)^j \left( {\frac{{L_r }}{{L_{r - 1} }}} \right)^j L_{r + j} }  = ( - 1)^k \left( {\frac{{L_r }}{{L_{r - 1} }}} \right)^k F_{k + 1} L_r
\end{equation}
and
\begin{equation}
\sum_{j = 0}^k {( - 1)^j \left( {\frac{{P_r }}{{P_{r - 1} }}} \right)^j P_{r  + j} }  = ( - 1)^k \left( {\frac{{P_r }}{{P_{r - 1} }}} \right)^k P_{k + 1} P_r\,.
\end{equation}
\begin{thm}\label{thm.ybopnqn}
For non-negative integer $k$, integers $m$ and $r$ and any integer $n$ for which $w_n\ne 0$, the following identity holds:
\[
q^{n-r} eu_{r - 1} \sum_{j = 0}^k {\frac{{u_{m - (n + 1) - kr + rj} }}{{(w_n /w_{n - r} )^j }}}  = \frac{{w_m w_{n - r} }}{{(w_n /w_{n - r} )^k }} - w_n w_{m - (k + 1)r}\,. 
\]

\end{thm}
In particular,
\begin{equation}\label{eq.ndpr9xm}
q^{n-r} eu_{r - 1} \sum_{j = 0}^k {\frac{{u_{rj} }}{{(w_n /w_{n - r} )^j }}}  = \frac{{w_{n + kr + 1} w_{n - r} }}{{(w_n /w_{n - r} )^k }} - w_n w_{n - r + 1}\,.
\end{equation}
\begin{proof}
Write identity~\eqref{eq.vx6b1t4} as
\[
w_m  = \frac{{w_n }}{{w_{n - r} }}w_{m - r}  + q^{n - r} \frac{{eu_{r - 1} }}{{w_{n - r} }}u_{m - n - 1}\,. 
\]
Identify $x=w_n/w_{n-r}$, $y=q^{n-r}eu_{r-1}/w_{n-r}$, $\alpha=r$ and $\beta=n+1$ and use these in Lemma~\ref{lem.u4bqbkc}. 
\end{proof}
Results for the $\{G_m\}$ and $\{P_m\}$ sequences emanating from identity~\eqref{eq.ndpr9xm} are the following:
\begin{equation}
( - 1)^{n - r} (G_0 G_1  + G_0^2  - G_1^2 )F_r \sum_{j = 0}^k {\frac{{F_{rj + 1} }}{{(G_n /G_{n - r} )^j }}}  = \frac{{G_{n + kr + 1} G_{n - r} }}{{(G_n /G_{n - r} )^k }} - G_n G_{n - r + 1}
\end{equation}
and
\begin{equation}
( - 1)^{n - r - 1} P_r \sum_{j = 0}^k {\frac{{P_{rj + 1} }}{{(P_n /P_{n - r} )^j }}}  = \frac{{P_{n + kr + 1} P_{n - r} }}{{(P_n /P_{n - r} )^k }} - P_n P_{n - r + 1}\,.
\end{equation}
\begin{lemma}\label{lem.s9jfs7n}
Let $\{X_m\}$ be any arbitrary sequence, where $X_m$, $m\in\Z$, satisfies a second order recurrence relation $X_m=xX_{m-\alpha}+yX_{m-\beta}$, where $x$ and $y$ are arbitrary non-vanishing complex functions, not dependent on $m$, and $\alpha$ and $\beta$ are integers. Then,
\begin{equation}\label{eq.mxyb9zk}
y\sum_{j = 0}^k {\frac{{X_{m - k\alpha  - \beta  + \alpha j} }}{{x^j }}}  = \frac{{X_m }}{{x^k }} - xX_{m - (k + 1)\alpha }\,, 
\end{equation}
\begin{equation}\label{eq.cgldajj}
x\sum_{j = 0}^k {\frac{{X_{m - k\beta  - \alpha  + \beta j} }}{{y^j }}}  = \frac{{X_m }}{{y^k }} - yX_{m - (k + 1)\beta }\,,
\end{equation}
\begin{equation}\label{eq.n2n4ec3}
\sum_{j = 0}^k { \frac{X_{m - (\beta - \alpha)k + \alpha + (\beta - \alpha)j}}{(-y/x)^j} }  = \frac{x X_m}{(-y/x)^k}  + y X_{m - (k + 1)(\beta - \alpha)}
\end{equation}
and
\begin{equation}\label{eq.c522g7v}
\sum_{j = 0}^k { \frac{X_{m - (\alpha - \beta)k + \beta + (\alpha - \beta)j}}{(-x/y)^j} }  = \frac{y X_m}{(-x/y)^k}  + x X_{m - (k + 1)(\alpha - \beta)}\,.
\end{equation}
for $k$ a non-negative integer.
\end{lemma}
In particular,
\begin{equation}
y\sum_{j = 0}^k {\frac{{X_{\alpha j} }}{{x^j }}}  = \frac{{X_{k\alpha+\beta} }}{{x^k }} - xX_{\beta-\alpha }\,, 
\end{equation}
\begin{equation}
x\sum_{j = 0}^k {\frac{{X_{\beta j} }}{{y^j }}}  = \frac{{X_{k\beta+\alpha} }}{{y^k }} - yX_{\alpha -\beta }\,,
\end{equation}
\begin{equation}
\sum_{j = 0}^k { \frac{X_{(\beta  - \alpha )j}}{(-y/x)^j} }  =  \frac{xX_{(\beta  - \alpha )k - \alpha }}{(-y/x)^k}  + yX_{ - \beta }
\end{equation}
and
\begin{equation}
\sum_{j = 0}^k { \frac{X_{(\alpha  - \beta )j}}{(-x/y)^j} }  =  \frac{yX_{(\alpha  - \beta )k - \beta }}{(-x/y)^k}  + xX_{ - \alpha }\,.
\end{equation}
\begin{proof}
Identity~\eqref{eq.mxyb9zk} is a direct consequence of Lemma~\ref{lem.u4bqbkc} with $Y_m=X_m$. Identity~\eqref{eq.cgldajj} is obtained from the symmetry of the recurrence relation by interchanging $x$ and $y$ and $\alpha$ and $\beta$. Identity~\eqref{eq.n2n4ec3} is obtained from~\eqref{eq.mxyb9zk} by re-arranging the recurrence relation and comparing coefficients. Identity~\eqref{eq.c522g7v} is obtained from identity~\eqref{eq.n2n4ec3} by interchanging $x$ and $y$ and $\alpha$ and $\beta$. 
\end{proof}
Note that the identities~\eqref{eq.mxyb9zk}, \eqref{eq.cgldajj} and~\eqref{eq.n2n4ec3} can be writen in the following equivalent forms:
\begin{equation}\label{eq.awbhgnm}
y\sum_{j = 0}^k {x^j X_{m - \beta  - \alpha j} }  = X_m  - x^{k + 1} X_{m - (k + 1)\alpha }\,,
\end{equation}
\begin{equation}\label{eq.jjikwds}
x\sum_{j = 0}^k {y^j X_{m - \alpha  - \beta j} }  = X_m  - y^{k + 1} X_{m - (k + 1)\beta }\,,
\end{equation}
\begin{equation}
\sum_{j = 0}^k {\frac{{X_{m + \alpha  - (\beta  - \alpha )j} }}{{( - x/y)^j }}}  = xX_m  + \frac{y}{{( - x/y)^k }}X_{m - (k + 1)(\beta  - \alpha )}
\end{equation}
and
\begin{equation}
\sum_{j = 0}^k {\frac{{X_{m + \beta  - (\alpha  - \beta )j} }}{{( - y/x)^j }}}  = yX_m  + \frac{x}{{( - y/x)^k }}X_{m - (k + 1)(\alpha  - \beta )}\,.
\end{equation}
In particular,
\begin{equation}
y\sum_{j = 0}^k {x^j X_{- \alpha j} }  = X_\beta  - x^{k + 1} X_{\beta - (k + 1)\alpha }\,,
\end{equation}
\begin{equation}
x\sum_{j = 0}^k {y^j X_{- \beta j} }  = X_\alpha  - y^{k + 1} X_{\alpha - (k + 1)\beta }\,,
\end{equation}
\begin{equation}
\sum_{j = 0}^k {\frac{{X_{(\alpha  - \beta )j} }}{{( - x/y)^j }}}  = xX_{ - \alpha }  + \frac{y}{{( - x/y)^k }}X_{k(\alpha  - \beta ) - \beta }
\end{equation}
and
\begin{equation}
\sum_{j = 0}^k {\frac{{X_{(\beta  - \alpha )j} }}{{( - y/x)^j }}}  = yX_{ - \beta }  + \frac{x}{{( - y/x)^k }}X_{k(\beta  - \alpha ) - \alpha }\,.
\end{equation}
\begin{thm}\label{thm.m4yo6a8}
For non-negative integer $k$ and any integer $m$, the following identities hold:
\begin{equation}\label{eq.vybd467}
  qu_r^ku_{r - 1} \sum_{j = 0}^k {\frac{{w_{m - kr - r - 1 + rj} }}{{u_r^j }}}  = u_r^{k+1} w_{m - kr - r} - w_m,\quad r\in\Z,\quad r\ne-1\,,
\end{equation}
\begin{equation}\label{eq.vwqo0w9}
u_{r-1} \sum_{j = 0}^k {\frac{{w_{m - kr - r + 1 + rj} }}{{(-qu_{r - 2})^j }}}  = \frac{{w_m }}{(-qu_{r - 2})^k} + qu_{r - 2} w_{m - (k + 1)r},\quad r\in\Z,\quad r\ne1\,, 
\end{equation}
and
\begin{equation}\label{eq.utwljqu}
\sum_{j = 0}^k {\frac{{w_{m - k + r + j} }}{{(qu_{r - 1} /u_r )^j }}}  = \frac{{u_rw_m }}{{(qu_{r - 1} /u_r )^k }} - qu_{r - 1}w_{m - k - 1},\quad r\in\Z,\quad r\ne0\,.
\end{equation}

\end{thm}
In particular,
\begin{equation}
qu_r^ku_{r - 1} \sum_{j = 0}^k {\frac{{w_{rj} }}{{u_r^j }}}  =bu_r^{k+1} - w_{kr + r + 1}\,,
\end{equation}
\begin{equation}
u_{r-1} \sum_{j = 0}^k {\frac{{w_{rj} }}{(-qu_{r - 2})^j}}  =  \frac{{w_{kr + r - 1} }}{(-qu_{r - 2})^k} + (ap - b)u_{r - 2}
\end{equation}
and
\begin{equation}\label{eq.btkvoap}
\sum_{j = 0}^k {\frac{{w_j }}{{(qu_{r - 1} /u_r )^j }}}  = \frac{{u_rw_{k - r} }}{{(qu_{r - 1} /u_r )^k }} - \frac{1}{{q^r }}\frac{{au_{r + 1}  - bu_r }}{{au_{r + 1}  + (b - pa)u_r }}u_{r - 1}w_{r + 1}\,.
\end{equation}
\begin{proof}
Write the relation~\eqref{eq.fuxige6} as \mbox{$w_m  = u_r w_{m - r}  - qu_{r - 1} w_{m - r - 1}$}, identify $X=w$, $x=u_r$, $y=-qu_{r-1}$, $\alpha=r$ and $\beta=r+1$ and use these in Lemma~\ref{lem.s9jfs7n}.

\end{proof}
Explicit examples from identity~\eqref{eq.btkvoap} include
\begin{equation}
\sum_{j = 0}^k {( - 1)^j \frac{{G_j }}{{(F_r /F_{r + 1} )^j }}}  = ( - 1)^k \frac{{F_{r + 1} }}{{(F_r /F_{r + 1} )^k }}G_{k - r}  - ( - 1)^r \frac{{F_{r + 2} G_0  - F_{r + 1} G_1 }}{{F_{r + 2} G_0  + F_{r + 1} (G_1  - G_0 )}}F_rG_{r + 1}\,,
\end{equation}
\begin{equation}
\sum_{j = 0}^k {( - 1)^j \frac{{P_j }}{{(P_r /P_{r + 1} )^j }}}  = ( - 1)^k \frac{{P_{r + 1} P_{k - r} }}{{(P_r /P_{r + 1} )^k }} + ( - 1)^r P_r P_{r + 1}
\end{equation}
and
\begin{equation}
\sum_{j = 0}^k {\frac{{( - 1)^j }}{{2^j }}\frac{{J_j }}{{(J_r /J_{r + 1} )^j }}}  = \frac{{( - 1)^k }}{{2^k }}\frac{{J_{r + 1} J_{k - r} }}{{(J_r /J_{r + 1} )^k }} + \frac{{( - 1)^r }}{{2^r }}J_r J_{r + 1}\,.
\end{equation}
\begin{thm}\label{thm.yng8u8b}
For non-negative integer $k$ and all integers $r$ and $m$, the following identities hold
\begin{equation}\label{eq.u5k6v3w}
\sum_{j = 0}^k {\frac{{w_{m - kr + r + rj} }}{{(q^r /v_r )^j }}}  = \frac{{v_r w_m }}{{(q^r /v_r )^k }} - q^r w_{m - (k + 1)r}\,,
\end{equation}
\begin{equation}\label{eq.x6yh3ef}
v_r^kq^r \sum_{j = 0}^k {\frac{{w_{m - r + rj} }}{{v_r^j }}} = v_r^{k+1} w_m   - w_{m + (k + 1)r}
\end{equation}
and
\begin{equation}\label{eq.is4vgui}
v_r \sum_{j = 0}^k {\frac{{w_{m - 2kr - r + 2rj} }}{{( - q^r )^j }} = \frac{{w_m }}{{( - q^r )^k }} + q^r w_{m - (k + 1)2r} }\,.
\end{equation}

\end{thm}
In particular,
\begin{equation}
\sum_{j = 0}^k {\frac{{w_{rj} }}{{(q^r /v_r )^j }}}  = \frac{{v_r w_{kr - r} }}{{(q^r /v_r )^k }} - \frac{1}{{q^r }}\frac{{au_{2r}  - bu_{2r - 1} }}{{au_{2r}  + (b - pa)u_{2r - 1} }}w_{2r}\,,
\end{equation}
\begin{equation}
v_r^kq^r \sum_{j = 0}^k {\frac{{w_{rj} }}{{v_r^j }}} = v_r^{k+1} w_r   - w_{(k + 2)r}
\end{equation}
and
\begin{equation}
v_r \sum_{j = 0}^k {\frac{{w_{2rj} }}{{( - q^r )^j }} = \frac{{w_{(2k + 1)r} }}{{( - q^r )^k }} + \frac{{au_r  - bu_{r - 1} }}{{au_r  + (b - pa)u_{r - 1} }}w_r }\,.
\end{equation}
\begin{proof}
Write identity~\eqref{eq.w7u7hr6} as \mbox{$w_m=(1/v_r)w_{m+r}+(q^r/v_r)w_{m-r}$}. Identify $X=w$, $x=1/v_r$, $y=q^r/v_{r}$, $\alpha=-r$ and $\beta=r$ and use these in Lemma~\ref{lem.s9jfs7n}, identities~\eqref{eq.cgldajj}, \eqref{eq.n2n4ec3} and~\eqref{eq.awbhgnm}.
\end{proof}
\begin{lemma}\label{lem.i84yg3s}
Let $\{X_m\}$ be any arbitrary sequence. Let $X_m$, $m\in\Z$, satisfy a second order recurrence relation $X_m=xX_{m-\alpha}+yX_{m-\beta}$, where $x$ and $y$ are non-vanishing complex functions, not dependent on $m$, and $\alpha$ and $\beta$ are integers. Then,
\begin{equation}\label{eq.nrzg4pd}
\sum_{j = 0}^k {\binom kj\left( {\frac{x}{y}} \right)^j X_{m - k\beta  + (\beta  - \alpha )j} }  = \frac{{X_m }}{{y^k }}\,,
\end{equation}
\begin{equation}\label{eq.h6kcv7w}
\sum_{j = 0}^k {\binom kj\frac{{X_{m + (\alpha - \beta)k + \beta j} }}{{( - y )^j }}}  = \left( { - \frac{x}{y}} \right)^k X_m
\end{equation}
and
\begin{equation}\label{eq.fnwrzi3}
\sum_{j = 0}^k {\binom kj\frac{{X_{m + (\beta - \alpha)k + \alpha j} }}{{( - x )^j }}}  = \left( { - \frac{y}{x}} \right)^k X_m\,,
\end{equation}
for $k$ a non-negative integer.

\end{lemma}
In particular,
\begin{equation}
\sum_{j = 0}^k {\binom kj\left( {\frac{x}{y}} \right)^j X_{(\beta  - \alpha) j} }  = \frac{{X_{k\beta} }}{{y^k }}\,,
\end{equation}
\begin{equation}
\sum_{j = 0}^k {\binom kj\frac{{X_{\beta j} }}{{( - y)^j }}}  = \left( { - \frac{x}{y}} \right)^k X_{(\beta  - \alpha )k}
\end{equation}
and
\begin{equation}
\sum_{j = 0}^k {\binom kj\frac{{X_{\alpha j} }}{{( - x)^j }}}  = \left( { - \frac{y}{x}} \right)^k X_{(\alpha  - \beta )k}\,.
\end{equation}
\begin{proof}
Only identity~\eqref{eq.nrzg4pd} needs to be proved, since identity~\eqref{eq.h6kcv7w} is obtained from~\eqref{eq.nrzg4pd} by re-arranging the recurrence relation while identity~\eqref{eq.fnwrzi3} is obtained from identity~\eqref{eq.h6kcv7w} by interchanging $\alpha$ and $\beta$ and $x$ and $y$. To prove~\eqref{eq.nrzg4pd}, we apply mathematical induction on $k$. Obviously, the theorem is true for $k=0$. We assume that it is true for $k=n$ a positive integer. The induction hypothesis is
\[
\mathcal{P}_n :\left(f(n)  = \frac{{X_m }}{{y^n }}\right)\,;
\]
where
\[
f(n)=\sum_{j = 0}^n {\binom kj\left( {\frac{x}{y}} \right)^j X_{m - n\beta  + (\beta  - \alpha )j} }.
\]
We want to prove that $\mathcal{P}_n\Rightarrow\mathcal{P}_{n+1}$. We proceed,
\[
f(n + 1) = \sum_{j = 0}^{n + 1} {\binom{n+1}j\left( {\frac{x}{y}} \right)^j X_{m - n\beta  - \beta  + (\beta  - \alpha )j} }\qquad\qquad\qquad\qquad\qquad\qquad\qquad\qquad\qquad 
\]
\[
\mbox{ (since $\binom{n+1}j=\binom nj+\binom n{j-1}$})
\]
\[
 = \sum_{j = 0}^{n + 1} {\binom nj\left( {\frac{x}{y}} \right)^j X_{m - n\beta  - \beta  + (\beta  - \alpha )j} }  + \sum_{j = 0}^{n + 1} {\binom n{j-1}\left( {\frac{x}{y}} \right)^j X_{m - n\beta  - \beta  + (\beta  - \alpha )j} } 
\]
\[
 = \sum_{j = 0}^{n + 1} {\binom nj\left( {\frac{x}{y}} \right)^j X_{m - n\beta  - \beta  + (\beta  - \alpha )j} }  + \sum_{j = 1}^{n + 1} {\binom n{j-1}\left( {\frac{x}{y}} \right)^j X_{m - n\beta  - \beta  + (\beta  - \alpha )j} } 
\]
\[
\quad = \sum_{j = 0}^n {\binom nj\left( {\frac{x}{y}} \right)^j X_{m - n\beta  - \beta  + (\beta  - \alpha )j} }  + \frac{x}{y}\sum_{j = 0}^n {\binom n{j}\left( {\frac{x}{y}} \right)^j X_{m - n\beta  - \beta  + (\beta  - \alpha )(j + 1)} } 
\]
\[
 = \sum_{j = 0}^n {\binom nj\left( {\frac{x}{y}} \right)^j \left( {X_{m - n\beta  - \beta  + (\beta  - \alpha )j}  + \frac{x}{y}X_{m - n\beta  - \beta  + (\beta  - \alpha )(j + 1)} } \right)}\qquad\qquad 
\]
\[
 = \frac{1}{y}\sum_{j = 0}^n {\binom nj\left( {\frac{x}{y}} \right)^j \left( {xX_{m - n\beta  - \beta  + (\beta  - \alpha )(j + 1)}  + yX_{m - n\beta  - \beta  + (\beta  - \alpha )j} } \right)}\qquad\quad
\]
(since $xX_{m - n\beta  - \beta  + (\beta  - \alpha )(j + 1)}  + yX_{m - n\beta  - \beta  + (\beta  - \alpha )j}  = X_{m - n\beta  + (\beta  - \alpha )j} $)
\[
 = \frac{1}{y}\sum_{j = 0}^n {\binom nj\left( {\frac{x}{y}} \right)^j X_{m - n\beta  + (\beta  - \alpha )j} }\qquad\qquad\qquad\qquad\qquad\qquad\qquad\quad 
\]
\[
\quad\qquad\qquad = \frac{1}{y}\frac{{X_m }}{{y^n }}\quad\mbox{(by the induction hypothesis)}\,.\qquad\qquad\qquad\qquad\qquad\qquad\qquad\quad
\]
Thus,
\[
f(n + 1) = \frac{{X_m }}{{y^{n + 1} }}\,,
\]
so that
\[
\mathcal{P}_{n+1}:\left(f(n+1)  = \frac{{X_m }}{{y^{n+1} }}\right)\,; 
\]
i.~e.~$\mathcal{P}_n\Rightarrow\mathcal{P}_{n+1}$ and the induction is complete.
\end{proof}
Note that the identity of Lemma~\ref{lem.i84yg3s} can also be written as
\begin{equation}
\sum_{j = 0}^k {\binom kj\left( {\frac{y}{x}} \right)^j X_{m - k\alpha  + (\alpha  - \beta )j} }  = \frac{{X_m }}{{x^k }}\,,
\end{equation}
with the particular case
\begin{equation}
\sum_{j = 0}^k {\binom kj\left( {\frac{y}{x}} \right)^j X_{(\alpha  - \beta) j} }  = \frac{{X_{k\alpha} }}{{x^k }}\,.
\end{equation}
\begin{thm}\label{thm.g1ihfq5}
For non-negative integer $k$ and any integer $m$, the following identities hold:
\begin{equation}\label{eq.f9x35z3}
(-qu_{r - 1})^k\sum_{j = 0}^k {\binom kj\left( {-\frac{{u_r }}{{qu_{r - 1} }}} \right)^j w_{m - k(r + 1) + j} }  = w_m,\quad r\in\Z,\quad r\ne0\,,
\end{equation}
\begin{equation}\label{eq.r5w2cg1}
\sum_{j = 0}^k {\binom kj \frac{w_{m - k + rj}}{(qu_{r - 2})^j} }  = \left( {\frac{{u_{r-1} }}{{qu_{r - 2} }}} \right)^k w_m, \quad r\in\Z,\quad r\ne1,
\end{equation}
and
\begin{equation}\label{eq.fxtzfk3}
\sum_{j = 0}^k {(-1)^j\binom kj \frac{w_{m + k + rj}}{u_r^j} }  = \left( {\frac{{qu_{r - 1} }}{{u_r }}} \right)^k w_m,\quad r\in\Z,\quad r\ne-1\,. 
\end{equation}

\end{thm}
In particular,
\begin{equation}
(-qu_{r - 1})^k\sum_{j = 0}^k {\binom kj \left( {-\frac{{u_r }}{{qu_{r - 1} }}} \right)^j w_j }  = w_{k(r + 1)}\,,
\end{equation}
\begin{equation}\label{eq.wbtbfxw}
\sum_{j = 0}^k {\binom kj\frac{w_{rj}}{(qu_{r - 2})^j} }  = \left( {\frac{{u_{r-1} }}{{qu_{r - 2} }}} \right)^k w_k
\end{equation}
and
\begin{equation}
\sum_{j = 0}^k {(-1)^j\binom kj \frac{w_{rj}}{u_r^j} }  = \left( {\frac{{u_{r - 1} }}{{u_r }}} \right)^k \frac{{au_k  - bu_{k - 1} }}{{au_k  + (b - pa)u_{k - 1} }}w_k\,. 
\end{equation}

\begin{proof}
Use, in Lemma~\ref{lem.i84yg3s}, the $x$, $y$, $\alpha$ and $\beta$ found in the proof of Theorem~\ref{thm.m4yo6a8}. 
\end{proof}
We have the following specific examples from identity~\eqref{eq.wbtbfxw}:
\begin{equation}
\sum_{j = 0}^k {( - 1)^j \binom kj\frac{{G_{rj} }}{{F_{r - 1}^j }}}  = ( - 1)^k \left( {\frac{{F_r }}{{F_{r - 1} }}} \right)^k G_k\,,
\end{equation}
\begin{equation}
\sum_{j = 0}^k {( - 1)^j \binom kj\frac{{P_{rj} }}{{P_{r - 1}^j }}}  = ( - 1)^k \left( {\frac{{P_r }}{{P_{r - 1} }}} \right)^k P_k
\end{equation}
and
\begin{equation}
\sum_{j = 0}^k {\frac{{( - 1)^j }}{{2^j }}\binom kj\frac{{J_{rj} }}{{J_{r - 1}^j }}}  = \frac{{( - 1)^k }}{{2^k }}\left( {\frac{{J_r }}{{J_{r - 1} }}} \right)^k J_k\,.
\end{equation}
Note that identity~\eqref{eq.f9x35z3} is a generalization of identity~(48) of Vajda~\cite{vajda}, the latter being the evaluation of the former at $r=1,q=-1$.
\begin{thm}\label{thm.peyb26i}
For non-negative integer $k$ and all integers $m$ and $r$, the following identities hold:
\begin{equation}\label{eq.e6qnu1m}
\sum_{j = 0}^k {\binom kj\frac{{w_{m - kr + 2rj} }}{{q^{rj} }}}  = \left( {\frac{{v_r }}{{q^r }}} \right)^k w_m\,,
\end{equation}
\begin{equation}\label{eq.k130vx8}
\sum_{j = 0}^k {\binom kj\left( { - \frac{{v_r }}{{q^r }}} \right)^j w_{m - 2kr + rj} }  = \frac{{w_m }}{{( - q^r )^k }}
\end{equation}
and
\begin{equation}\label{eq.d00yx5i}
\sum_{j = 0}^k {(-1)^j\binom kj\frac{{w_{m + kr + rj} }}{v_r^j }}  = (-1)^k\frac{{q^{rk} w_m }}{{v_r^k }}\,.
\end{equation}
\end{thm}
In particular,
\begin{equation}
\sum_{j = 0}^k {\binom kj\frac{{w_{2rj} }}{{q^{rj} }}}  = \left( {\frac{{v_r }}{{q^r }}} \right)^k w_{rk}\,,
\end{equation}
\begin{equation}
\sum_{j = 0}^k {\binom kj\left( { - \frac{{v_r }}{{q^r }}} \right)^j w_{rj} }  = \frac{{w_{2kr} }}{{( - q^r )^k }}
\end{equation}
and
\begin{equation}\label{eq.xf5dcmx}
\sum_{j = 0}^k {(-1)^j\binom kj\frac{{w_{rj} }}{v_r^j}}  = (-1)^k\left( {\frac{{au_{kr}  - bu_{kr - 1} }}{{au_{kr}  + (b - pa)u_{kr - 1} }}} \right)\frac{{w_{kr} }}{{ v_r^k}}\,.
\end{equation}
\begin{proof}
Use, in Lemma~\ref{lem.i84yg3s}, the $x$, $y$, $\alpha$ and $\beta$ found in the proof of Theorem~\ref{thm.yng8u8b}.
\end{proof}
Setting $p=1=-q$ in identity~\eqref{eq.e6qnu1m}, we have
\begin{equation}
\sum_{j = 0}^k {( - 1)^{rj} \binom kjG_{m - kr + 2rj} }  = ( - 1)^{rk} L_r^k G_m\,.
\end{equation}
Identity~\eqref{eq.xf5dcmx} at $p=1=-q$ gives
\begin{equation}
\sum_{j = 0}^k {( - 1)^j \binom kj\frac{{G_{rj} }}{{L_r^j }}}  = ( - 1)^k \frac{{F_{kr + 1} G_0  - F_{kr} G_1 }}{{F_{kr + 1} G_0  + F_{kr} (G_1  - G_0 )}}\frac{{G_{kr} }}{{L_r^k }}\,.
\end{equation}
\section{Acknowledgement}
I would like to thank the anonymous referee whose comments helped to improve this paper.


\begin{thebibliography}{99}
\bibitem{horadam61} A.~F.~Horadam, \emph{A generalized Fibonacci sequence}, The American Mathematical Monthly, \textbf{68.5} (1961), 455--459.
\bibitem{horadam65} A.~F.~Horadam, \emph{Basic properties of a certain generalized sequence of numbers}, The Fibonacci Quarterly, \textbf{3.3} (1965), 161--176.
\bibitem{stanica03} P.~Stanica, \emph{Generating functions, weighted and non-weighted sums for powers of secondorder recurrence sequences}, The Fibonacci Quarterly, \textbf{41.4} (2003), 321--333.
\bibitem{vajda} S.~Vajda, \emph{Fibonacci and Lucas Numbers, and the Golden Section: Theory and Applications}, Dover Press, (2008).


\end{thebibliography}
\end{document}